\magnification=1200
\centerline{\bf Lengths of geodesics between two points on a Riemannian manifold.}
\bigskip
\centerline{\bf Alexander Nabutovsky and Regina Rotman}
\bigskip
\noindent
{\bf Abstract.} Let $x$ and $y$ be two (not necessarily distinct) points
on a closed Riemannian manifold $M^n$. According to a well-known theorem
by J.P. Serre there exist infinitely many geodesics between $x$ and $y$.
It is obvious that the length of a shortest of these geodesics cannot exceed 
the diameter of the manifold. But what can be said about the lengths of
the other geodesics? We conjecture that for every $k$ there are $k$ distinct
geodesics of length $\le k\ diam(M^n)$. This conjecture is evidently
true for round spheres and is not difficult to prove for all closed
Riemannian manifolds with non-trivial torsion-free fundamental groups. 
In this paper we announce two further results in the
direction of this conjecture: Our first result is
that there always exists a second geodesic between $x$ and $y$ of length not
exceeding $2n\ diam(M^n)$.  Our second result is 
that if $n=2$ and $M^2$ is diffeomorphic to $S^2$, then for
every $k$ every
pair of points of $M^2$ can be connected by $k$ distinct geodesics of length
less than or equal to
$(4k^2-2k-1)diam(M^2)$.
\bigskip\noindent
{\bf 1. Introduction} A classical theorem of J.P. Serre ([S])
established the existence of infinitely many geodesics between any two
points on a closed Riemannian manifold. Yet 
almost nothing is known about lengths of these geodesics.
Of course, it is almost a tautology that the length of a shortest of
these geodesics does not exceed the diameter of the manifold.
But what can be said about the lengths of a second shortest geodesic,
the third shortest geodesic, etc. ? We think that the following
conjecture could be true:
\medskip\noindent
{\bf Conjecture 1.} Let $x$ and $y$ be two points in a closed Riemannian
manifold $M^n$. Then for every $k=1,2,3,\ldots$ there exist $k$ distinct
geodesics starting at $x$ and ending at $y$ of length not exceeding
$k\ diam(M^n)$.
\bigskip
Note that this conjecture is obvious for the round spheres. On the other
hand it is not difficult to prove that:
\bigskip\noindent
{\bf Proposition 2.} Assume that the fundamental group of $M^n$
is non-trivial and torsion-free.
Then Conjecture 1 is true for $M^n$.
\bigskip\noindent
{\bf Proof.} Consider a non-contractible loop $\gamma$ based at $x$ of
length $\le 2diam(M^n)$, realizing a non-trivial
element of $\pi_1(M^n)$ of infinite order. (According to [Gr]
there exists a basis for $\pi_1(M^n)$ made of loops of length
$\leq 2diam(M^n)$ based at $x$.)
Let $p$ denote
the midpoint of $\gamma$. Let $\tau$ denote a shortest geodesic
that goes from $x$ to $y$ and $\sigma$ denote a shortest geodesic
that goes from $y$ to $p$. Finally, let $\gamma_1$ and $\gamma_2$ denote
two distinct geodesics that go from $x$ to $p$ along $\gamma$.
Without
any loss of generality we can assume that the loop $\gamma_1*\sigma^{-1}*
\tau^{-1}$ is not contractible. (Here and below for every
path $\rho$ $\rho^{-1}$ denotes
$\rho$ traversed in the opposite direction.)
Now $\tau$, $\gamma_1*\sigma^{-1}$, $\gamma^{-1}*\tau$, $\gamma*\gamma_1*
\sigma^{-1}$, $\ldots$ , $\gamma^{-i}*\tau$,
$\gamma^i*\gamma_1*\sigma^{-1}$, $\ldots$ are pairwise non-path homotopic
paths connecting $x$ and $y$. The lengths of the first $k$ of them do not
exceed $k\ diam(M^n)$. We obtain the desired geodesics minimizing the length
in the corresponding path homotopy classses. QED.
\bigskip
The main purpose of this paper is to announce the following two results:
\bigskip\noindent
{\bf Theorem 3.} Let $M$ be a closed
Riemannian manifold diffeomorphic
to $S^2$. Let $x$ and $y$ be two (not necessarily distinct) points of $M$.
Then for every $k=1,2,\ldots$ there exist $k$ geodesics starting
at $x$ and ending at $y$ of length not exceeding
$(4k^2-2k-1)diam(M^2)$.
If $x=y$, then this upper bound can be replaced by
a better bound
$(4k^2-6k+2)diam(M^2)$. (In this case $k$ distinct geodesics are $k$ geodesic loops
based at $x$; $k-1$ of them are non-trivial).
\bigskip\noindent
{\bf Theorem 4.} Let $M^n$ be a closed Riemannian manifold of dimension $n$.
Let $q$ be the minimal $i$ such that $\pi_i(M^n)\not= 0$. Let $x,y$ be
any two points of $M^n$. Then there exist two distinct geodesics
between $x$ and $y$ such that their lengths do not exceed
$2q\ diam(M^n)(\le 2n\ diam(M^n))$.
\bigskip\noindent
{\bf Remark.}
Assume that $x=y$ in Theorem 4. Then the shortest closed geodesic
between $x$ and $y$ is trivial. Therefore in this case the statement of Theorem 4
asserts that for every $x\in M^n$ there exists a non-trivial geodesic
loop in $M^n$ based at $x$ of length $\leq 2q\ diam(M^n)\leq 2n\ diam(M^n)$.
This assertion had been proven by the second author in [R].
\medskip
A complete proof of Theorem 3 will appear in [NR1], and a complete
proof of Theorem 4 will appear in [NR2].
\bigskip\noindent
{\bf 2. A sketch of the proof of Theorem 4.}
\medskip\noindent
If $x=y$, then the shortest geodesic between $x$ and $y$ is the trivial
geodesic, and the theorem asserts that the length of a second shortest
geodesic loop based at an arbitrary point $x$ of $M^n$ does not
exceed $2qd\leq 2nd$. Here and below $d$ denotes the diameter of $M^n$.
This result has been recently proven by the second
author in [R]. The proof in the general case when $x\not= y$
generalizes the proof in the case when $x=y$.
\par
We can reinterprete the main construction in [R] as follows: Assume
that there are no geodesic loops of length $\leq L$ based at $x$,
and therefore the space $\Omega_L^x$ of all loops of length $\leq L$
based at $x$ is contractible. This assumption can be used in order
to perform the following inductive construction:
Let $z$ be any point of $M^n$, and
$C_i^z$ denotes the space of
all $i$-tuples of paths between $x$ and $z$ of length
$\leq {L\over 2(i-1)}$ for some $i=2,\ldots$. Then, using
an induction with respect to $i$  one can construct a
continuous map from $C_i^z$ into the space of continuous maps
from the $i$-dimensional disc $D^i$ into $M^n$. Each
$(i+1)$-tuple of paths between $x$ and $z$ becomes a collection
of meridians of the image of $\partial D^i$. So, the $i$-disc
fills the $i$-tuple of paths. This construction desribed in [R]
has a number of good properties. In particular, consider
$D^i$ as a regular $i$-dimensional simplex. Then, if one applies
this construction to a subset that is made of some $(i-1)$ paths among
the $i$ paths forming an $i$-tuple, then the resulting singular
$(i-1)$-simplex will be a $(i-1)$-face
of the singular $i$-simplex corresponding
to the $i$-tuple. This
construction can be generalized in order to fill not only $i$-tuples
of paths between $x$ and $z$ but 1-dimensional complexes with
vertices $x,z_1,\ldots,z_i$ pairwise connected by paths so that
the lengths of all paths between $z_i$ and $z_j$ are very small. (In other
words the point $z$ ``splits" into $i$ very close points connected
by very short paths.)
\par
Now the upper bound $2qd$ for the length of the shortest geodesic loop
based at $x$ can be proven by contradiction.  Assume that there are no
such loops, so that the construction outlined above can be applied. Take
a non-contractible map of a $q$-sphere into $M^n$. Extend it to a map
of a $(q+1)$-disc thereby obtaining a contradiction as follows: Triangulate
$D^{q+1}$ as a cone over a very fine triangulation of $S^q$. Map the
center of $D^{q+1}$ into $x$, and all new 1-simplices into minimizing
geodesics that start at $x$. Now extend the constructed map of the
1-skeleton of every $(q+1)$-simplex $\sigma$ of the triangulation of
$D^{q+1}$
to a map of $\sigma$ by applying the construction above.
\par
It turns out that in order to generalize this proof to the situation
when $x\not= y$ one only needs to change the base of the inductive
construction
of maps of $C_i^z$ into $Map(D^i, M^n)$, namely the construction
of the map of $C_2^z$ into $Map(D^2, M^n)$. In [R] this filling of
digons was based on the following easy lemma:  Let $\gamma_1,\gamma_2$ be
two paths between $x$ and $z$  of length $\leq l$. Assume that there are
no non-trivial
geodesic loops based at $x$ of length $\leq 2l$. Then there exists
a path homotopy between $\gamma_1$ and $\gamma_2$ that passes through paths
between $x$ and $z$ of length $\leq 3l$.
\par
So, our goal can be achieved by replacing this lemma by the following lemma:
\medskip\noindent
{\bf Lemma 5.} Let $l$ be a positive number,
$x,y,z$ points in a Riemannian manifold $M^n$,
and $\gamma_1$, $\gamma_2$ two paths of length $\leq l$
between $x$ and $z$ in $M^n$. Assume that
the distance between $y$ and $z$ does not exceed $l$.
Finally, assume that there exists exactly one geodesic
of length $\leq 2l$ between $x$ and $y$. Then $\gamma_1$ and $\gamma_2$ can
be connected by a path homotopy that passes only through paths between
$x$ and $z$ of length $\leq 3l$.
\medskip
\medskip\noindent
{\bf 3. A sketch of the proof of Theorem 3.}
\medskip\noindent
Let us start from the following definition:\medskip\noindent
{\bf Definition 6.} Let $M^n$ be a
Riemannian manifold diffeomorphic to $S^n$, and $L$
a positive number. Let $f:S^n\longrightarrow M^n$ be a map of
a non-zero degree such that the length of the image of every meridian of
$S^n$ does not exceed $L$. Then we call $f$ a {$L$-controlled
vertical sweep-out of $M^n$}. 
\medskip
The following proposition easily follows from a  well-known
geometric description of generators of homology groups
of the loop space of $S^n$ (cf. [Schw]):
\medskip\noindent
{\bf Proposition 7.}  Assume that $M^n$ is diffeomorphic to $S^n$
and admits an $L$-controlled vertical sweep-out $f$ for some $L$.
Let $x$ be the image of
the South pole of $S^n$ under $f$, and $y$ be an arbitrary point of $S^n$.
Then for every positive integer $k$ there exist
$k$ distinct geodesics between $x$ and $y$ of length
$\leq 2(k-1)L+diam(M^n)$. If, in addition, $x=y$ then this upper bound can be replaced by
a better bound $2(k-1)L$.
\medskip
In particular, the proposition implies that if $M^n$ is diffeomorphic
to $S^n$ then the length of the $k$th geodesic between $x$ and
$y$ is bounded by a linear function of $k$. The proof of this well-known
fact for all simply connected closed Riemannian manifolds can be found
in [Schw].  This proposition has the following immediate corollary:
If we were able to derive an upper
bound for $L$ in an optimal vertical sweep-out  in terms of $diam(M^n)$
(and possibly $n$), then we would get an upper bound for the length
of the $k$th geodesic between $x$ and $y$ of the form $c(n)k\ diam (M^n)$.
However, such an upper bound for $L$ cannot be, in general, true even if
$n=2$ and $M^n$ is diffeomorphic to $S^2$. (One can construct a family
of counterexamples using the example of M. Katz and S. Frankel [KF] of
a family of
metric 2-discs with diameter $D$ and the length of the boundary $l$ such
that in order to contract the boundary
one must first increase its length to $C(l+D)$, where $C$ can be an arbitrary
large constant.)
\par
We are going to explain the ideas of the proof only in the case when $M^2$ is a real
analytic Riemannian manifold. We refer the reader to the complete paper [NR1] for
the proof in the smooth case. 
\par
First, we are going to present a sketch of the proof in the case, when
$x=y$. We also assume that there exists a point $z\in M^2$ such that
$z$ is one of the most distant points for $x$, and $x$ is one of the most
distant points for $z$.
These assumptions enable one to prove the theorem in a more
geometric and transparent way than in the general case. We 
present here a rather detailed outline of the main geometric ideas
of the proof in this simpler case
in hope that this outline will facilitate the
understanding of the general proof in [NR1]. We will then
sketch how the proof in the simpler case can be modified to obtain
a proof of the theorem in general.
\par
Recall that the diameter of $M^2$ is denoted for brevity
by $d$.
In order to prove Theorem 3 we first attempt to
construct a $3d$-controlled vertical sweep-out of $M^2$ mapping the South pole
into $x$. We will see that our attempt
can only be prevented by a non-trivial geodesic
loop based at
$x$ of length $\leq 2d$. Then we will attempt to
construct a $5d$-controlled
vertical sweep-out so that our attempt can only be thwarted by a second
``short" geodesic loop based at $x$, and so on.
The worst situation appears when
the controlled sweep-out appears after
$k-2$ attempts that will be blocked by $k-2$ short non-trivial
geodesic loops based at $x$.
Then $L$ is proportional to $kd$, leading to a bound for the length of the
$k$th geodesic loop based at $x$ that is quadratic in $k$.
(In the worst case scenario the first $k-1$ geodesic loops based at $x$
provided by Proposition 7
can turn out to be $k-1$ short geodesic loops that were already obtained.)
Here are some details:
\par
First, we are going to use the assumption that $x$ is one of the most distant 
points for $z$, and $z$ is one of the most distant points for $x$.
According to an old observation of
\par\noindent
M. Berger in this case one can connect 
$x$ and $z$ by minimizing geodesics
so that all angles of geodesic digons formed by any two neighboring
geodesics $\gamma_i$ and $\gamma_{i+1}$
 do not exceed $\pi$. Consider
all minimizing geodesics between $x$ and $z$.
They split the manifold into digonal domains with angles $\leq\pi$.
When we apply the Birkhoff curve shortening process fixing the
endpoint $x$ to any
such digon the resulting homotopy takes place inside the domain it
bounds. As the
result we contract this digon either to $x$ or to a geodesic loop based at
$x$ and
contained in the closure of the domain bounded by this digon.
\par
If there are no geodesic loops 
of length $\leq 2d$ based at $x$,
then one can contract all such digons to $x$ as loops based
at $x$. Then one can use these homotopies to 
construct path homotopies between $\gamma_i$ and $\gamma_{i+1}$ passing
via curves between $x$ and $z$ of length $\leq 3d$. (See [R]
for details.) These path homotopies
together provide a $3d$-controlled vertical sweep-out of our 2-sphere.
\par
If there is no such a sweep-out, then the application of the
curve-shortening process to one such digon $D$ with vertices $x$, $z$
ends at a non-trivial geodesic loop $\gamma$ of length $\leq 2d$
based at $x$ inside the considered digonal domain.
Denote the midpoint of $\gamma$ by $m$.
If two segments of $\gamma$ between $x$ and $m$ can be connected by
a path homotopy without increase of length then $\gamma$ does not present
an obstruction to the attempted construction of the controlled
vertical sweep-out of the sphere. Moreover,
if there is such a path homotopy where lengths of paths are bounded
by some number $C$, then there exists a $(C+2d)$-controlled vertical
sweep-out. In order to construct such a controlled path homotopy
consider the part of the cut-locus of $x$ inside
the domain bounded by $\gamma$. (W.l.o.g. we can assume that $\gamma$ is
not a closed geodesic, since in this case we will obtain the desired
geodesic loops as multiples of this closed geodesic. By the domain bounded
by $\gamma$ we mean one of two connected components of the complement
of $\gamma$ that has an angle at $x$ less than $\pi$.) It is easy to see
that this part of the cut-locus is a non-empty finite tree, where there
exist exactly two minimizing geodesics connecting $x$ with any point on an edge of this
tree, and at least three minimizing geodesics connecting $x$ with any of its
vertices.
\par
If $m$ is not a vertex of this tree, it must be on an edge. Slide
it along the edge to the nearest vertex $v$ of the cut-locus inside the
domain bounded by $\gamma$. Correspondingly, we can find a homotopy of $\gamma$ 
to a digon $\gamma_1$ along digons of length $\leq 2d$ based at $x$
by connecting $x$ with the points on the edge with two geodesics.
All minimizing geodesics from $x$ to $v$ split a domain bounded by $\gamma_1$
into finitely many smaller domains with digonal boundaries
with angles at $x$ less than $\pi$.
All of these domains but possibly one have angles at $y$ not exceeding
$\pi$. Now we can apply the curve shortening process to boundaries 
of all these domains. For every digonal domain with angles $\leq\pi$ we
find ourselves in the situation considered above:
either its boundary can be contracted to $x$ by the curve-shortening process
, and we can eliminate it by contracting one of its edges to the other 
without significantly increasing its length, or we get stuck at
 a geodesic loop $\delta$, which is strictly inside of the domain bounded
by $x$, and so is a new geodesic loop. In this case we repeat the argument
for $\delta$ instead of $\gamma$, and so on. Every time when a new geodesic loop
appears, $L$ in the $L$-controlled vertical sweep-out
that we hope to construct can only increase by no more than $2d$.
If we can get rid of all digonal domains with angles $\leq\pi$ without
encountering a new geodesic loop, then $L$ does not increase.
\par
The only digonal domain,
 $\Delta$, where the angle at $v$ might be greater than $\pi$
seems to constitute a problem for us: When we apply the curve-shortening
process to its boundary the resulting geodesic
loop need not be inside this domain.
As the result the curve-shortening process can end not at a new geodesic
loop, but at a previously constructed geodesic loop, e.g. $\gamma$.
Let us call such a digonal domain {\it fat}.
So, we proceed in a slightly different manner. First, we observe that
since there are no geodesics between $x$ and $v$ inside $\Delta$, there is
an edge of the cut-locus of $x$ passing through $v$ or ending at $v$
inside $\Delta$ such that one can connect $x$ with points of this edge
in the closure of $\Delta$ by two continuously varying geodesics.
We will use these geodesics to slide $v$ and the boundary of $\Delta$ to
a digon between $x$ and another vertex $w$ of the cut-locus of $x$
strictly inside $\Delta$. Then we consider all minimizing geodesics between
$x$ and $w$ and
iterate the construction with $w$ instead of
$v$. In the absence of new obstructing geodesic loops the process of
constructing an $L$-controlled vertical sweep-out will end in finitely
many steps since the cut-locus of $x$ has finitely many edges.
\par
This argument is sufficient to take care of fat digons, and to
complete the proof in the case, when $x=y$ and $x$ is one
of the most distant from $z$ points. But we are going to
``improve" this argument having in mind the
general case. We are going to make the following observation
about another possible way to handle fat digons discussed in the
previous paragraph:
\par
Observe that the process of contracting digons and splitting digonal
domains at vertices of the cut locus can be represented by means of
a finite tree.
Let us call this tree {\it a filling tree}.
The initial digon is a root vertex of the tree.
Assume that our attempts to construct an $L$-controlled vertical
sweep-out are blocked by several obstructing geodesic loops before they
finally succeed as outlined above. Then $L\leq 3d+2d\lambda$, where
$\lambda$ denotes the maximal number of obstructing geodesic loops
that we can encounter along a path
from the root to one of the
vertices of degree one (=leaves) of the filling tree.
This observation enables us
to deal with the fat digonal domains in a somewhat different manner: We can
attempt to contract their boundaries using the Birkhoff curve-shortening
process. The cases when the fat digon contracts to a point, or a new
geodesic loop appears as the outcome of the Birkhoff curve-shortening
process can be treated exactly in exactly the same way as the case, when the
digon is not fat.
The trouble was that the result can be a geodesic loop $\omega$ that
has already arizen as an obstructing geodesic higher on the path from the
root of the tree that describes the filling process. But in this
case we can just connect $\omega$ with the fat digon by a homotopy passing
through loops of length not exceeding the length of the fat digon. (Observe
that this length does not exceed $2d$.)
The fat digon is inside the domain bounded by $\omega$, so the homotopy
between $\omega$ and the fat digon enables us to eliminate the whole
part of the filling tree between the point corresponding to $\omega$ and
the point corresponding to the fat digon. Of course, the value of $L$
does not increase as the result of homotoping $\omega$ to the fat digon.
Now we can consider the part of the cut locus of $x$ inside the domain
bounded by the fat digon (with the angle at $x$ less than $\pi$), slide
the fat digon to the nearest vertex, etc. Note that if we modify
the filling tree in accordance with the just explained idea, then
all obstructing geodesic loops encountered along any path will be distinct.
Therefore we will either obtain $k-1$ distinct short non-trivial geodesic
loops along some path (wich together with the trivial loop constitute
a set of $k$ geodesic loops we want to find), or obtain an $L$-controlled
vertical $L$ sweep-out of $M^2$ with  $L=(2k-1)d$. Now Proposition 7 implies the desired estimate
in the case, when $x=y$ and $x$ is one of the most distant points
for one of its most distant points.
\par
Let us very briefly outline what we need to change in this proof
to make it work in the general case when $x\not= y$
(or $x=y$ but $x$ is not one of the most distant for $z$ points of $M^2$).
We need to learn how
to contract loops based at $x$ using a curve shortening process for
paths between $x$ and $y$. Here is how it can be done: Fix a minimizing geodesic $\rho$ from $x$ to $y$.
Now proceed by means of the following path homotopies:
$$\gamma\longrightarrow\gamma*\rho*\rho^{-1}\longrightarrow \tau*\rho^{-1},$$
where $\tau$ is a geodesic between $x$ and $y$ obtained from $\gamma*\rho$ as the result of the application
of a curve shortening process on the space of paths from $x$ to $y$.
If $\tau=\rho$, then we can cancel $\rho*\rho^{-1}$ over itself, and we are done.
Otherwise, we will call $\tau$ an {\it obstructing geodesic}. If there are at most $k-1$ geodesics
between $x$ and $y$, then at most $k-2$ of them can arise as obstructing geodesics for various loops
based at $x$. Now the following observation is critical for our purposes: If two loops $\gamma_1$ and
$\gamma_2$ based at $x$ have the same obstructing geodesic, then they can be connected by a homotopy that
passes through loops of length $\leq\max\{length(\gamma_1),length(\gamma_2)\}+2dist(x,y)\leq\max\{length(\gamma_1),length(\gamma_2)\}+2d\leq 4d.$
Now it seems that proceeding as in the case, when $x=y$, that was described above, we either obtain $k$ distinct
obstructing geodesics between $x$ and $y$, or an $L$-controlled vertical sweep-out of $M^2$ with $L\leq 2(k-2)d+3d+2dist(x,y)=(2k-1)d+2dist(x,y)\leq (2k+1)d$. 
\par
Yet here we encounter a new technical problem. Previously, we used
the fact that the digonal angles formed by the neighbouring minimizing
geodesics between $x$ and $z$ had angles $\leq\pi$
to ensure that all path homotopies between these
minimizing geodesics sweep $M^2$
with degree one. (All path homotopies took places within domains bounded
by the considered digons.)
Now our homotopies are all over the place, and we do not have any
obvious way to ensure that all arcs constructed during our path homotopies
sweep out $M^2$ with multiplicity one.
\par
In order to circumvent this difficulty we are going to proceed as follows:
Our original strategy was to connect $x$ with one of the most distant point
from $x$ by minimizing geodesics, and to construct a vertical $L$-controlled
sweep-out of $M^2$ by contracting digons formed by neighboring pairs of
these minimizing geodesics. Now we are going to
abandon this plan in favor of the following
strategy (used in several previous papers of the authors): Consider a 
diffeomorphism $f:S^2\longrightarrow M^2$ and try to extend it to
a map of $D^3$ into $M^2$. Such an extension is obviously impossible.
Endow $S^2$ with a very fine triangulation, and triangulate $D^3$ as a
cone over the triangulation of $S^2$ with one new vertex, $p$. Map $p$
into $x$, map all new one dimensional simplices (``radii" of $D^3$) into
minimizing geodesics from $x$ to the images of the endpoints of the radii
under $f$. Before extending the map to new two-dimensional simplices, observe
that their boundaries where mapped into two minimizing geodesics emanating
from $x$ and a very short geodesic connecting the endpoints of these
minimizing geodesics. For all practical purposes these boundaries can
be be regarded as geodesic digons in $M^2$. We contract them
in the same way as we
described above (using the cut locus of $x$ to
contract loops). More precisely, we consider a domain on $M^2$ bounded
by the digon and the part of the cut locus of $x$ inside this domain.
For every vertex of the considered part of the cut locus (which is a tree)
we consider all minimizing geodesics from $x$ to
this vertex. These minimizing geodesics split the domain inside the
considered digon into smaller domains also bounded by digons, which are nested inside each other.
If we can find a path homotopy from one side of a digon  to the other
with a controlled increase of length we can ``eliminate" a domain
bounded by this digon. It is not difficult to see that it is sufficient to know
how to contract the boundaries of digons regarded as loop based at $x$ via loops based at $x$
of controlled length in order to find the required path homotopy.
We are, however, prevented from contracting these loops in an appropriate way
by obstructing geodesics. Yet, for every digon we can find a homotopy connecting it with
the loop that is the boundary of the furthest digon ``down the tree" with the same obstructing geodesic,
and the length of loops during this homotopy increases by no more than $2d$.
Therefore we find ourselves in the
situation, where each obstructing
geodesic influences the lengths of paths in only one path homotopy.
As above, this leads to an inductive process of contraction of loops based at
$x$ that provides a desired bound for the length of loops during the resulting homotopy.
This homotopy yields a desired filling of the considered 2-simplex.
\par
As the result, we either obtain $k$
distinct short obstructing geodesics between $x$ and $y$ ,
or construct a path homotopy
of one geodesic segment in the boundary of the digon to the other
that passes through short curves. If there are less than $k$ short obstructing
geodesics, we will construct an extension to the 2-skeleton of $D^3$.
Consider now the
boundaries of 3-simplices of $D^3$. The image of
one of four triangles in any such
boundary is very small, and for all practical purposes can be treated
as a point. For at least one of these boundary 2-spheres the constructed
map of this sphere in $M^2$ is not contractible. Because
of its construction we obtain the desired $L$-controlled vertical sweep-out
of $M^2$. Now an application of Proposition 7 completes the proof of the
theorem.
\bigskip\noindent
{\bf Acknowledgemets:} 
The work of both authors had been partially
supported by their NSF grants and their NSERC Discovery grants.
The work of Regina Rotman had been also
partially supported by NSERC University Faculty Award.
We would like to thank Dima Burago for his suggestions that
helped to improve the exposition.
\bigskip\noindent
{\bf References:}
\par\noindent
[Gr] M. Gromov, ``Metric structures for Riemannian and
 non-Riemannian spaces",\par\noindent Birkhauser, 1999.
\par\noindent
[KF] M. Katz, S. Frankel, ``The Morse landscape of a Riemannian disc",
Annales de l'Inst. Fourier, 43(2)(1993), 503-507.
\par\noindent
[NR1] A. Nabutovsky, R. Rotman, ``The length of geodesics on a two-dimensional sphere", preprint.
\par\noindent
[NR2] A. Nabutovsky, R. Rotman, ``The length of a
second shortest geodesic", preprint.
\par\noindent
[R] R. Rotman, ``The length of a shortest geodesic loop at a point",
submitted for publication, available at
http://comet.lehman.cuny.edu/sormani/others/rotman.html
\par\noindent
[S] J.P. Serre, ``Homologie singuli\`ere des espaces fibr\'es. Applications",
Ann. Math., 54(1951), 425-505.
\par\noindent
[Schw] A.S. Schwarz, ``Geodesic arcs on Riemannian manifolds", Uspekhi Math.
Nauk (translated from Russian as ``Russian Math. Surveys"),
13(6)(1958), 181-184.
\bigskip\noindent
Department of Mathematics, University of Toronto, Toronto, Ontario,
M5S2E4, CANADA and Department of Mathematics, McAllister Bldg.,
The Pennsylvania State University, University Park, PA 16802, USA.
\bye